\documentclass[a4paper,12pt,reqno]{amsart}
\usepackage{amssymb}
\usepackage[dvips]{graphicx}
\usepackage{hyperref}

\input epsf.tex
\newdimen\xsize
\newdimen\oldbaselineskip
\newdimen\oldlineskiplimit
\def\restorelineskip{\baselineskip=\oldbaselineskip%
\lineskiplimit=\oldlineskiplimit}
\def\putm[#1][#2]#3{
\hbox{\vbox to 0pt{\parindent=0pt%
\vskip#2\xsize\hbox to0pt{\hskip#1\xsize $#3$\hss}\vss}}}%
\long\def\Line#1{\hbox to \hsize{#1}}
\def\putt[#1][#2]#3{
\vbox to 0pt{\noindent\hskip#1\xsize\lower#2\xsize%
\vtop{\restorelineskip#3}\vss}}

\makeatletter
\def\xbig[#1]#2{{\hbox{$\m@th\left#2\vbox to#1\xsize{}%
\right.\n@space$}}}
\makeatother
\def\xlar[#1]#2{%
\smash{\mathop{ \hbox to #1\xsize{\leftarrowfill}}\limits^{#2}}}
\def\xrar[#1]#2{%
\smash{\mathop{ \hbox to #1\xsize{\rightarrowfill}}\limits^{#2}}}
\def\xline[#1]{\hbox to #1\xsize{\leaders\hrule\hfill}}

\DeclareFontFamily{U}{rsf}{\skewchar\font'177}%
\DeclareFontShape{U}{rsf}{m}{n}{<-6>rsfs5<6-8>rsfs7<8->rsfs10}{}%
\DeclareFontShape{U}{rsf}{b}{n}{<-6>rsfs5<6-8>rsfs7<8->rsfs10}{}%
\DeclareMathAlphabet\RSFS{U}{rsf}{m}{n}
\SetMathAlphabet\RSFS{bold}{U}{rsf}{b}{n}
\def\sf#1{{\mathsf{#1}}}

\def\slsf{\slshape \sffamily }

\def\msmall#1{\mathchoice{\hbox{\small$\displaystyle {#1}$}}{#1}{#1}{#1}}

\def\cc{{\mathbb C}}

\def\rr{{\mathbb R}}
\def\sss{{\mathbb S}}

\def\pp{{\mathbb P}}

\def\adyn{\sf{1}}

\def\dim{\sf{dim}\,}

\def\dist{\sf{dist}\,}

\def\span{\sf{span}}

\def\lim{\mathop{\sf{lim}}}

\def\max{\sf{max}}

\def\supp{\sf{supp}\,}

\def\v{{\mathrm{v}}}
\def\w{{\mathrm{w}}}

\def\eps{\varepsilon}

\def\<{\langle}\let\la=\<
\def\>{\rangle}\let\ra=\>

\def\d{\partial}

\def\ddef{\mathrel{{=}\raise0.3pt\hbox{:}}}
\def\deff{\mathrel{\raise0.3pt\hbox{\rm:}{=}}}

\def\fraction#1/#2{\mathchoice{{\msmall{ #1\over#2}}}%
{{ #1\over #2 }}{{#1/#2}}{{#1/#2}}}

\def\le{\leqslant}

\def\emptyset{\varnothing}

\def\longpoints{\leaders\hbox to 0.5em{\hss.\hss}\hfill \hskip0pt}
\def\stateskip{\smallskip}
\def\state#1. {\stateskip\noindent{\bf#1. }} 
\def\statep#1. {\stateskip\noindent{\bf#1 }} 
\def\proof{\state Proof. \2}

\def\Chi{\raise 2pt\hbox{$\chi$}}

\def\ie{\hskip1pt plus1pt{\sl i.e.\/,\ \hskip1pt plus1pt}}

\def\sli{{\sl i)} } 
\def\slii{{\sl i$\!$i)} } 
\def\sliii{{\sl i$\!$i$\!$i)} }
\def\sliv{{\sl i$\!$v)} }

\def\Chi{\raise 2pt\hbox{$\chi$}}

\let\phI=\phi\let\phi=\varphi\let\varphi=\phI
\let\cal=\mathcal

\def\calb{{\cal B}}
\def\calc{{\cal C}}
\def\cald{{\cal D}}

\def\calg{{\cal G}}

\def\calk{{\cal K}}

\def\calo{{\cal O}}

\def\calu{{\cal U}}
\def\calv{{\cal V}}
\def\calw{{\cal W}}
\def\calx{{\cal X}}
\def\caly{{\cal Y}}

\def\eps{\varepsilon}

\def\d{\partial}

\def\1{{1\mkern-5mu{\rom l}}}

\def\ge{\geqslant}

\def\fraction#1/#2{\mathchoice{{\msmall{ #1\over#2}}}%
{{ #1\over #2 }}{{#1/#2}}{{#1/#2}}}

\def\le{\leqslant}

\def\emptyset{\varnothing}
\newcommand{\2}{\thinspace}

\def\qed{\ \ \hfill\hbox to .1pt{}\hfill\hbox to .1pt{}\hfill $\square$\par}

\def\comment#1\endcomment{}

 \belowdisplayskip=\abovedisplayskip

\def\lineeqqno(#1){\hfill\llap{\vbox to 10pt%
{\vss\begin{align} \eqqno(#1)\end{align}\vss}}\vskip1pt}

\advance\headheight 1.2pt

\def\ShowwLLabel#1{}

\def\thechpt{\Roman{chpt}}

\def\newchapt[#1]#2{\newpage%
\refstepcounter{chpt}\setcounter{subsection}{0}%
\setcounter{thm}{0}\setcounter{defi}{0}%
\setcounter{rema}{0}\setcounter{exrc}{0}%
\renewcommand{\thesubsection}{\thechpt.\arabic{subsection}}%
\section*{\begin{center}\huge \bf Chapter \thechpt\\
#2 \end{center}}\label{#1}%
\ \smallskip%
\markboth{Chapter \thechpt}{#2}%
}
\def\newsect[#1]#2{\refstepcounter{section}\setcounter{equation}{0}%
\renewcommand{\thesubsection}{\arabic{section}.\arabic{subsection}}%
\section*{\arabic{section}.
#2}\vspace{-20pt}\label{#1}\vspace{20pt}%
\markboth{Section \arabic{section}}{#2}}

\def\newlect[#1]#2{\refstepcounter{section}%
\renewcommand{\thesubsection}{\arabic{section}.\arabic{subsection}}%
\section*{Lecture \arabic{section}\\
#2}\label{#1}%
\markboth{Lecture \arabic{section}}{#2}}

\def\newprg[#1]#2{\refstepcounter{subsection}%
\subsection*{{\thesubsection.\ #2}} \label{#1}%
}

\def\newappx[#1]#2{%
\refstepcounter{appx}\setcounter{section}{0}%
\renewcommand{\thesubsection}{A\arabic{appx}.\arabic{subsection}}%
\section*{Appendix \arabic{appx}\\ #2}
\label{#1}%
\markboth{Appendix A\arabic{appx}}{#2}
}

\newtheorem{thm}{Theorem}[section]
   \def\newthm#1{\begin{thm}\label{#1}}

\newtheorem{nnthm}{Theorem}
   \def\newthm#1{\begin{nnthm}\label{#1}}

\newtheorem{lem}{Lemma}[section]
   \def\newlemma#1{\begin{lem} \label{#1}}

\newtheorem{prop}{Proposition}[section]
   \def\newprop#1{\begin{prop}\label{#1}}

\newtheorem{nnprop}{Proposition}
   \def\newprop#1{\begin{nnprop}\label{#1}}

\newtheorem{corol}{Corollary}[section]
   \def\newcorol#1{\begin{corol} \label{#1}}

\newtheorem{nncorol}{Corollary}
   \def\newcorol#1{\begin{nncorol} \label{#1}}

\newtheorem{defi}{Definition}[section]
   \def\newdefi#1{\begin{defi} \label{#1}\rm }

\newtheorem{nndefi}{Definition}
   \def\newdefi#1{\begin{nndefi} \label{#1}\rm }

\newtheorem{exmp}{Example}[section]
   \def\newexmp#1{\begin{exmp} \label{#1}\rm }

\newtheorem{exrc}{Exercise}
   \def\newexrc#1{\begin{exrc} \label{#1}\rm }

\newtheorem{rema}{Remark}[section]
   \def\newrema#1{\begin{rema} \label{#1}\rm }

   \def\newrema#1{\begin{rema} \label{#1}\rm }

\def\eqqno(#1){\label{(#1)}}
\def\eqqref(#1){(\ref{(#1)})}

\usepackage[active]{srcltx}   

\def\el2{\sf{L^2}}
\def\el{\sf{l}}

\title{Loop  Spaces as Hilbert-Hartogs Manifolds}
\author{M. Anakkar, S. Ivashkovich}
\address{*
Universit\'e de Lille-1, UFR de Math\'ematiques, 59655 Villeneuve
d'Ascq, France.}
\email{ivachkov@math.univ-lille1.fr , mohammed.anakkar@univ-lille.fr}

\thanks{* Supported by the Labex Cempi ANR-11-LABX-0007-01}
\date{\today}
\date{\today}

\begin{document}

\begin{abstract}
We prove that generalized loop spaces of Hartogs manifolds are Hilbert-Hartogs.
We prove also that Hilbert-Hartogs manifolds possess a better extension properties 
that it is postulated in their definition. Finally, we give a list of examples
of Hilbert-Hartogs manifolds.
\end{abstract}

\maketitle

\setcounter{tocdepth}{1}
\tableofcontents

\newsect[INT]{Introduction and statement of results}

\newprg[INT.h-h]{Hilbert-Hartogs manifolds}

Let $X$ be a finite dimensional complex manifold, $q$ a natural number. We say that $X$ is $q$-Hartogs 
if every holomorphic mapping of a $q$-concave Hartogs figure 
\begin{equation}
\eqqno(q-1-fig)
H_q^1(r):= \left(\Delta^q\times \Delta (r)\right)\cup \left(A^q_{1-r,1} \times\Delta \right)
\subset \cc^{q+1}
\end{equation}
with values in $X$ extends to the unit polydisc $\Delta^{q+1}$. Here $\Delta^q(r)$ is a polydisk in 
$\cc^q$ of radius $r>0$ and $A^q_{r_1,r_2}=\Delta^q (r_2)
\setminus \bar \Delta^q(r_1)$ is a ring domain, $r_2>r_1$. Analogously to the finite dimensional 
case we say that a complex Hilbert manifold $\calx$ is a $q$-Hartogs Hilbert
manifold (or simply Hilbert-Hartogs when $q$ is clear from the context) if every holomorphic mapping 
$f:H_q^1(r) \to \calx$ 
of the $q$-concave Hartogs figure with values in $\calx$ extends to a holomorphic map 
$\tilde f:\Delta^{q+1}\to \calx$ of the unit polydisk to $\calx$. Hilbert manifolds in this paper
are modeled over $\el^2$. The latter stands for the space of the square integrable
sequences of complex numbers with its natural scalar product $(z,w) = \sum_kz_k\bar w_k$ and the norm 
$||z||^2 = \sum_k|z_k|^2$. In addition our manifolds, both finite and infinite dimensional, are supposed 
to be second countable. First of all let us remark that Hilbert-Hartogs manifolds possess much stronger 
extension properties than postulated in their definition.

\begin{nndefi}
A $q$-concave Hartogs figure in $\el^2$ is the following open set
\begin{equation}
\eqqno(hart-q-infty)
H_q^{\infty}(r):= \left(\Delta^q\times B^{\infty}(r)\right)\cup \left(A^q_{1-r ,1} \times
B^{\infty}\right),
\end{equation}
\noindent where $0< r <1$.
\end{nndefi}
Here $B^{\infty}(r)$ stands for the ball in $\el^2$ of radius $r>0$ centered at origin. When writing
$D\times \calb $, where $D$ is a subset in $\cc^q$ with coordinates $z=(z_1,...,z_q)$ and $\calb$ is a subset in 
$\el^2$ with coordinates $w= (w_1,....)$ we always mean the obvious subset in $\el^2$ with coordinates 
$(z,w) = (z_1,...,z_q,w_1,...)$, namely $D\times \calb = \{(z,w)\in \el^2:z\in D, w\in \calb \}$.

\begin{nnthm}
\label{q1-qinfty}
Let $\calx$ be a $q$-Hartogs Hilbert manifold. Then 
every holomorphic mapping $f:H_q^{\infty}(r)\to \calx$ extends to a
holomorphic mapping $\tilde f: \Delta^q\times B^{\infty}\to \calx$.
\end{nnthm}
As an immediate consequence  we obtain the following statement. A domain $\cald$ in a complex 
Hilbert manifold $\calx$ with the boundary of class $\calc^2$ is called $q$-pseudoconcave at its 
boundary point $p$ if the complex Hessian of a local defining $\cald$ function at $p$ has at least 
$q$ negative eigenvalues when restricted to $T^c_p\d \cald$. Here by saying that $\rho$ defines 
$\cald $ at $p\in \d \cald$ we mean that for some neighborhood $\calu$ of $p$ one has 
$\cald\cap \calu = \{x:\rho (x) <0\}$ (and of course $\nabla \rho|_{\d\cald}$ never vanishes).

\begin{nncorol}
\label{q-levi}
Let $\cald$ be a domain in a complex Hilbert manifold $\calx$ which is $q$-pseudocon\-cave at a 
boundary point $p$. Then every holomorphic map $f: \cald \to \caly$ to a $q$-Hartogs Hilbert manifold
$\caly$ extends holomorphically to a neighborhood of $p$.
\end{nncorol}

The proof of this corollary follows from Theorem \ref{q1-qinfty} by appropriately placing the Hartogs figure 
$H_{q}^{\infty}(r)$ near $p\in \cald$.

\newprg[INT.loop]{Loop spaces as Hilbert-Hartogs manifolds}
Our main goal in this paper is to show that one finds Hilbert-Hartogs manifolds more often than 
one could expect. In order to provide such examples we concentrate our attention in section 
\ref{LOOP} on loop spaces of complex manifolds. Fix a smooth real compact  manifold 
$S$, with or without boundary, and a finite dimensional complex manifold $X$. Then the manifold 
$\calw^{k,2}_{S,X}\deff W^{k,2}(S,X)$
of Sobolev $W^{k,2}$-maps carries a natural structure of a complex Hilbert manifold, see \cite{L1} or, 
section \ref{LOOP}. Here $k\ge n=\dim_{\rr}S$ to ensure that mappings from this space  are continuous. 
The manifold $\calw^{k,2}_{S,X}$ is usually called a {\slsf generalized loops space} of $X$. We prove the 
following statement.

\begin{nnthm}
\label{loop-hart}
A generalized loop space of a $q$-Hartogs complex manifold is a $q$-Hartogs Hilbert manifold.
\end{nnthm}
This provides us a lot of interesting examples of infinite dimensional Hilbert-Hartogs manifolds. For 
example let $S=\sss^1$ be the circle and $X$ a connected Riemann surface. Remark that $X$ is $1$-Hartogs if and 
only if $X$ is different from the Riemann sphere $\pp^1$ (see section \ref{LOOP} for details). Therefore all 
loop spaces $\calw^{1,2}_{\sss^1,X}$ 
for $X\not=\pp^1$ are Hilbert-Hartogs. Another example is obtained by lettting $G$ be a complex Lie group. 
Then by \cite{ASY} $G$ is Hartogs.
Therefore the loop space, which is classically denoted as $LG$, is also Hartogs ($\adyn$-Hartogs to be precise).
From the results of \cite{Iv1,Iv2,IS}, which determine Hartogs manifolds of complex dimensions $2$ and $3$ 
we get the following statement.

\begin{nncorol}
\label{sept-0}
Let $X$ be a compact complex manifold of dimension $2$ (resp. of dimension $3$).
Then: 

\smallskip\sli either $X$ contains a $2$-dimensional spherical shell (resp. a $3$-dimensional shell),
or a 

\quad rational curve (resp. a uniruled complex surface),

\smallskip\sliii or, every generalized loop space $\calw^{k,2}_{S,X}$ is $1$-Hartogs (resp. $2$-Hartogs).
\end{nncorol}

\noindent $X$ can contain both shells and rational curves.
Both these objects are obstructions to holomorphic extension. But if $X$ doesn't contain neither shells no rational
curves then $X$ is Hartogs ($1$-Hartogs if $\dim X=2$ and $2$-Hartogs if $\dim X=3$).
It might be interesting to think about $X$ from this Corollary as being (an unknown) surface of class $VII_0^{+}$ 
or as $\sss^6$ (provided it exists as a complex manifold), see subsection \ref{LOOP.almost-h} for more details and explanations.

\smallskip\noindent{\slsf Acknowledgment.} The Authors would like to express their gratitude to our colleague 
L\'ea Blanc-Centi for useful discussions around the results of this paper. We would like also to thank the Referee
of this paper who pointed to us a serious gap in its first version.

\newsect[HH]{Hilbert-Hartogs Manifolds}

\newprg[HH.proof]{$1$-complete neighborhoods}

In this subsection we recall and make a bit more precise the results from \cite{AZ}. Let $\calu$ 
be an open subset of $\el^2$.  Let furthermore  $u \in \mathcal{C}^2(\calu, \rr)$.

\begin{defi}
\label{psh-f}
We say that $u$ is {\slsf strictly plurisubharmonic} if for every $z\in \calu$
\[
\sum_{k,l=1}^{\infty}{\frac{\partial^2u}{\partial z_k \partial\bar{z}_l}(z)v_k\bar{v}_l} 
\geqslant c(z)||v||^2 \quad \forall v\in \el^2,
\]
with $c(z)$ strictly positive and continuous on $\calu$.
\end{defi}

Let $\calx$ be a complex Hilbert manifold modeled over $\el^2$ and let $\calv$  be an open subset 
in $\calx$. A function $u:\calv \to \rr$ is called strictly plurisubharmonic if for every coordinate 
chart $\phi_j:\calv_j\to \calu_j\subset \el^2$ the composition $u\circ \phi_j^{-1}$ is strictly 
plurisubharmonic in $\calu_j$.

\medskip 
By a closed coordinate ball we shall mean a closed subset $\bar\calb $ of $\calx$ such that there exists a 
coordinate chart $(\calv_j, \phi_j, \calu_j)$ with 

\medskip 1) $\bar B^{\infty} \subset \calu_j$;

\smallskip 2) $\bar\calb = \phi_j^{-1}(\bar B^{\infty})$, where $\bar B^{\infty}$
is the closed unit ball in $\el^2$.

\medskip In this situation $\calb\deff \phi_j^{-1}(B^{\infty})$ will be called an open coordinate ball, or
simply a coordinate ball. Here $B^{\infty}$ is the open unit ball in $\el^2$. Let $\calk$ be a {\slsf connected} 
compact in $\calx$.

\begin{defi}
\label{1-compl}
By a $1$-complete neighborhood of $\calk$ we understand an open set $\calv \supset \calk$ such that:

\smallskip\sli $\calv$ is contained in a finite union of open coordinate balls centered at points of 
$\calk$, \ie  
 
 \[
\calv \subset \bigcup_{j=1}^N\calv_j \quad\text{ with } \quad \calv_j = \phi_j^{-1}(B^{\infty})
\quad\text{ and }\quad \phi_j^{-1}(0) = k_j\in \calk;
 \]

\smallskip\slii $\calv$ possesses a strictly plurisubharmonic exhaustion function 
 $u:\calv\to [0, t_0)$, \ie 

\begin{itemize}
 \item for every $t<t_0$ one has that $\overline{u^{-1}\left([0,t)\right)}\subset \calv$,
 
 \item for all $0\le t_1<t_2<t_0$ one has that $\overline{u^{-1}\left([0,t_1)\right)}\subset
 u^{-1}\left([0,t_2)\right)$,
 
 \item $\calv = u^{-1}\left([0,t_0)\right)$ and $u^{-1}(0)$ is a single point in $\calk$.
\end{itemize}
\end{defi}

\noindent All closures are taken in the topology of $\calx$. We shall need the following result from \cite{AZ}.
By an analytic $q$-disk in a complex Hilbert manifold $\calx$ we understand an image $\phi (\bar D)$ of the 
closure of a bounded strictly pseudoconvex domain $D\subset \cc^q$ under a holomorphic imbedding $\phi :\bar D
\to \calx $, which is defined in a neighborhood of $\bar D$.

\begin{thm}
\label{a-z}
Let $\phi :\bar D\to \calx$ be an imbedded analytic $q$-disk in a second countable complex Hilbert
manifold modeled over $\el^2$. Then $\phi (\bar D)$ admits a fundamental system of $1$-complete 
neighborhoods.
\end{thm}

\noindent The proof of \ref{a-z} follows the main lines of the first part of Royden's proof from \cite{Ro}.
After trivializing the tubular neighborhood of $\phi (\bar D)$ to a sufficiently high degree one considers 
for every $0<\lambda <1$ a
function $\phi_{\lambda}(u,v) = \theta (u) + \lambda^{-1}||v||^2$. Here $u$ is a coordinate ``along''
$\phi (D)$ and $\theta$ is a strictly plurisubharmonic. exhaustion function of $\bar D = \{\theta (u) \le 1\}$, $v$ is 
a ``normal'' to $\phi (D)$ coordinate. The coordinates $u$ and $v$ are not holomorphic but only smooth
ones obtained by patching the local holomorphic coordinates. Using a sufficient degree of triviality of 
the tubular neighborhood one proves that these $\phi_{\lambda}$ are strictly plurisubharmonic. And then one defines $1$-complete 
neighborhoods of $\phi (\bar D)$ as follows: $\calo_{\lambda} = \{\phi_{\lambda}(u,v)<1+\lambda \}$.

\newprg[HH.inf-hart]{Infinite dimensional Hartogs figures}

Now we shall prove Theorem \ref{q1-qinfty} from the Introduction. For a unit vector $\v\in\el^2$ orthogonal to $\cc^q$ 
set $L_{\v} \deff \span \left\lbrace e_1,...,e_q,\v\right\rbrace$. Here we identify $\cc^q$ with the subspace 
$\el^2_q$ of $\el^2$ generated by $e_1,...,e_q$. Note that $L_{\v}\cap H_q^{\infty}(r)=H^1_q(r)$ and therefore given a 
holomorphic mapping $f:H_q^{\infty}(r) \to \calx$ its restriction to $L_{\v}\cap H_q^{\infty}(r)$ holomorphically 
extends to $L_{\v}\cap (\Delta^q\times B^{\infty})$. We conclude that for every line $<\v>\perp\cc^q$ the
restriction $f|_{L_{\v}}$ holomorphically extends to $L_{\v}\cap (\Delta^q\times B^{\infty})$, giving us an extension 
$\tilde f$ of $f$ to $\Delta^q\times B^{\infty}$. This extension is correctly defined because for
unit vectors $\v\not=\w$ orthogonal to $\cc^q$ the spaces $L_{\v}$ and $L_{\w}$ intersect only along $\cc^q\times \{0\}$. 

\smallskip Let us prove the continuity of $\tilde f$. Consider a sequence $(Z_n)_{n \geqslant 1}$ in 
$\Delta^q\times B^{\infty}$. Write it as $Z_n=(z^n , w^n)$ with $z^n\in \Delta^q$ and $w^n\in B^{\infty}$.
Suppose that $Z_n$ converges to $Z_0 = (z^0,w^0)\in \Delta^q\times B^{\infty}$.  Take 
$R$ such that $1-r<R<1$ with $||z^n||,||w^n||<R$ for all $n \in \mathbb{N}$, as well as $||z^0||,||w^0||<R$. 
Let $\phi_n: \bar{\Delta}^{q+1}\to \mathcal{X}$ be an analytic $(q+1)$-disk defined by $\phi_n(z,\eta)=
\tilde{f}(Rz,\eta w^n)$. Theorem \ref{a-z}   gives a $1$-complete neighborhood $V$ of the graph of 
$\phi_0 (z,\eta)=\tilde{f}(Rz, \eta w^0)$ over $\bar\Delta^{q+1}$. For $w^n$ close enough to $w^0$
the graph of $\phi_n$ over $L_{w_n} \cap H_q^\infty(r)$ is contained in $V$ because $L_{w_n} \cap H_q^\infty(r) 
\subset H_q^\infty(r)$ where $\tilde{f}$ is holomorphic. By the maximum principle for subharmonic functions
the graph of $\phi_n$ over the whole polydisk $\bar\Delta^{q+1}$ 
is contained in the neighborhood $V$ for $n$ big enough. This implies that $\phi_n$ converge uniformly to $\phi_0$ 
on $\bar\Delta^{q+1}$. Indeed, since $\phi_0$ is uniformly continuous we can find for a given $\eps >0$ a 
$\delta >0$ such that $d(\phi_0(x),\phi_0(y)) <\eps$ for all $x, y\in \bar\Delta^{q+1}$ such that 
$||x-y|| <\delta$. Here $d$ is a Riemannian metric 
compatible with the topology on $\calx$. Assume without loss of generality that $\delta <\eps$. For $n\gg 
N(\delta)$ there exists for every $x\in \bar\Delta^{q+1}$ an
 $y_n\in \bar\Delta^{q+1}$ such that $\dist \left((x,\phi_n(x)),(y_n,\phi_0(y_n)\right) <\delta$,
here $\dist$ is the product distance on $\Delta^{q+1}\times\calx$. But then $||x-y_n|| <\delta$ and therefore
$d(\phi_0(x),\phi_0(y_n))<\eps$. As the result we get that 
\[
\dist\left((x, \phi_n(x)), (x,\phi_0(x)\right)<
\dist\left((x, \phi_n(x)), (y_n,\phi_0(y_n)\right) +
\]
\[
+ \dist\left((y_n,\phi_0(y_n), (x,\phi_0(x)\right)
< \delta + \delta + \eps <3\eps .
\]
Therefore $\tilde f (z^0,w^0) = \phi_0\left(\frac{z^0}{R}, 1\right)=\lim_{n\to\infty}\phi_n\left(\frac{z^n}{R},1\right)
= \lim_{n\to\infty}\tilde f(z^n,w^n)$, \ie $\tilde{f}$ is continuous.

\smallskip  Now let us prove that this extension is G\^ateaux differentiable. Take again some $Z_0=(z^0,w^0)\in \Delta^q\times B^{\infty}$ and fix some direction $\v =(\v_1,\v_2)$ in $\cc^q\times\el^2$. Let $l= \{Z_0 + t\v:t\in\cc\}$ be the line through $Z_0$ in the direction $\v$. Set $L \deff \span\{e_1,...,e_q,\v_2\}$. Remark that $L$ contains $\cc^q$
and $l$. Note also that 
\begin{equation}
\eqqno(q-2-fig)
L\cap H_q^{\infty}(r) = H_q^2(r):= \left(\Delta^q\times B^2(r)\right)\cup \left(A^q_{1-r,1} \times B^2 \right)
\subset \cc^{q+2}
\end{equation}
is a Hartogs figure of bidimension $(q,2)$. The restriction $g\deff\tilde f|_{\Delta^q\times B^2}$ is continuous 
and holomorphic on $H_q^2(r)$. Holomorphicity of $g$ on the whole $\Delta^q\times B^2$  is a local question
and follows from the classical Hartogs separate analyticity theorem. Therefore $\tilde f$ is holomorphic on
$\Delta^q\times B^2$. In particular it is differentiable in the direction of $\v$ at $Z_0$, \ie is G\^ateaux differentiable. Every continuous G\^ateaux differentiable map is holomorphic, see Theorem 8.7 in \cite{Mu}, 
therefore Theorem \ref{q1-qinfty} is proved.

\smallskip\qed

\begin{rema}\rm
One can also consider the following version of an infinite dimensional Hartogs figure:

\begin{equation}
\eqqno(hart-inf-inf)
H_{\infty}^{\infty}(r) = \left(B^{\infty}\times B^{\infty}(r)\right)\cup
\left(\left(B^{\infty}\setminus \bar B^{\infty}(1-r)\right)\times B^{\infty}\right).
\end{equation}
Analogously to the proof of theorem above one can prove the following.
\begin{prop}
\label{infty-infty}
Holomorphic maps from $H_{\infty}^{\infty}(r)$ to $q$-Hartogs
Hilbert manifolds extend to $B^{\infty}\times B^{\infty}$
for any value of  $q \geqslant 1$. 
\end{prop}
Indeed, let
$f : H_\infty^\infty(r) \to \mathcal{X}$.
Choose vectors $e_1,...,e_{q-1}$ in $\el^2$ and let $\el^2_{q-1}$ denote their linear span. For any 
$v \in \left(\el^2_{q-1}\right)^{\perp}$, one has:
\begin{equation}
( \span\{e_1,...,e_{q-1},v\}\times \el^2) \cap H_\infty^\infty(r) = H_q^\infty(r).
\end{equation} 
Then, by Theorem \ref{q1-qinfty} $f$ extends along $L_v$ for all $v \in \left(\el^2_{q-1}\right)^{\perp}$ to
$\tilde{f}:B^\infty \times B^\infty \to \calx$.
The proof of continuity and Gateaux differentiability is the same as above.
\end{rema}

\newsect[LOOP]{Loop spaces of Hartogs manifolds are Hilbert-Hartogs}

\newprg[LOOP.loop]{Loop spaces of complex manifolds}

Fix a compact, connected, $n$-dimensional smooth real manifold $S$ with or without boundary and
let us following \cite{L1} describe the natural complex Hilbert structure on the Sobolev manifold $W^{k,2}(S, X)$ 
of $W^{k,2}$-maps of $S$ to a finite dimensional complex manifold $X$. To speak about Sobolev spaces
it is convenient to suppose that $X$ is smoothly imbedded to some $\rr^N$ by Whitney Imbedding Theorem, 
see ex. \cite{Mi}. If $X$ is not compact, we suppose that this imbedding is proper. For the following basic 
facts about Sobolev spaces we refer to \cite{T}.

\medskip \sli $f\in W^{k,2}(\rr^n)$ $\Longleftrightarrow$ $(1+||\xi
||)^k\hat f\in L^2 (\rr^n)$, where $\hat f$ is the Fourier
transform of $f$.

\smallskip Moreover this correspondence is an isometry by
the Plancherel identity. One defines then for any positive real number $s$ the Sobolev-Slobodetskii space
$W^s (\rr^n)$=$\{ f: (1+||\xi ||)^s\hat f\in L^2 (\rr^n)\}$.

\smallskip
\slii If $s\ge\frac{n}{2}+\alpha $ with $0<\alpha <1$ then
$W^s(\rr^n)\subset \calc^{\alpha}(\rr^n)$ and this inclusion is a
compact operator. In particular $W^{n,2}(\rr^n)\subset
\calc^{\frac{1}{2}}(\rr^n)$.

\sliii If $s>\frac{n}{2}+k$  for a positive integer $k$, then $W^s(\rr^n)\subset
\calc^{k}(\rr^n)$.

\sliv If $0\le s<\frac{n}{2}$ then $W^s (\rr^n)\subset
L^{\frac{2n}{n-2s}}(\rr^n)$.

\medskip From (\slii one easily derives that if $f,g\in W^{n,2} (\rr^n)$ then
$fg\in W^{n,2}(\rr^n)$. This enables one to define correctly a $W^{k,2}$-vector bundle over an 
$n$-dimensional real manifold provided $k\ge n$. By that we mean that the transition functions 
of the bundle are in $W^{k,2}$. Condition $k\ge n$ will be always assumed from now on. 
Take now $f\in W^{k,2} (S,X)$. Note that by (\slii such $f$ is H\"older continuous. Consider the 
pullback $f^*TX$ as a complex Sobolev bundle over $S$. A neighborhood $V_f$ of the zero section of
$f^*TX\to S$ is an open set of the complex Hilbert space $W^{k,2}(S,f^*TX)$ of Sobolev sections of 
the pullback bundle. This $V_f$ can be naturally identified with a neighborhood of $f$ in $W^{k,2}
(S,X)$, thus providing a structure of complex Hilbert manifold on
$W^{k,2} (S,X)$; see \cite{F} and Lemma 2.1 in \cite{L1} for more details on this construction. 

\smallskip Another way to understand this complex structure on $W^{k,2} (S,X)$
is to describe what are analytic disks in $W^{k,2} (S,X)$.

\begin{lem}
\label{lemp}
Let $D$ and $X$ be finite dimensional complex manifolds and let $S$ be
an $n$-dimensional compact real manifold with boundary. A mapping $F:D\times S\to X$
represents a holomorphic map from $D$ to $W^{k,2} (S,X)$ (denoted by the
same letter $F$) if and only if the following holds:

\sli for every $s\in S$ the map $F(\cdot,s):D\to X$ is holomorphic;

\slii for every $z\in D$ one has $F(z,\cdot)\in W^{k,2} (S,X)$ and this
correspondence $D\ni z\to F(z,\cdot)\in W^{k,2} (S,X)$ is continuous
with respect to the Sobolev topology on $W^{k,2}(S,X)$ (and the standard
topology on $D$).
\end{lem}

For the proof we refer to \cite{L1}. Now let us prove the Theorem \ref{loop-hart}
from the Introduction.

\newprg[LOOP.hart]{Proof of Theorem \ref{loop-hart}}
Let $F:H_q^1(r)\to W^{k,2} (S,X)$ be a holomorphic map. We represent it as a map
$F:H_q^1(r)\times S\to X$ enjoying properties (\sli and (\slii of Lemma \ref{lemp} 
above. From the fact that $X$ is $q$-Hartogs we get that for every fixed
$s\in S$ the mapping $F(\cdot,s):H_q^1(r)\to X$ extends holomorphically to $F(\cdot,s):
\Delta^{q+1}\to X$ and we get a mapping $F: \Delta^{q+1} \times S\to X$. We need to 
prove that for every fixed $z\in \Delta^{q+1}$ one has that $F(z,\cdot)\in W^{k,2} (S,X)$ 
and that this correspondence is continuous. 

\smallskip Fix some  $z_0 \in \Delta^{q+1}$ and some $s_0 \in S$. Take $R<1$ such that $z_0 \in 
\Delta^{q+1}_R$. Let $g_{s_0}$ be the map to the graph $\Gamma_{F(*,s_0)}$ of  $F(*,s_0)$, \ie 
$g_{s_0}$ is defined by
\begin{equation}
g_{s_0} : \Delta^{q+1}\ni z \longrightarrow  (z, F(z, s_0))\in \Delta^{q+1} \times X.
\end{equation}
By Royden's Lemma  for $R<1$ there exists a holomorphic embedding $G:\Delta^{q+1}_R\times \Delta^m_R 
\to \Delta^{q+1} \times X $ such that $G(*,0)=g_{s_0}$, where $m=\dim (X)$. Then $V=G(\Delta^{q+1}_R 
\times \Delta^m_R)$ contains $\Gamma_{F(*,s_0)}$ over $\Delta^{q+1}_R$. Since on $H_q^1(r)$ the map 
$z \mapsto g_{*}(z) \in W^{k,2}(S,X)$ is continuous, there exists an $\epsilon >0$ 
such that for $s \in B(s_0,\epsilon)$ the graph $\Gamma(*,s)$ over $H_q^1(r)$ is contained in $V$. Therefore 
by the Hartogs theorem for holomorphic functions the graph $\Gamma(*,s)$ over $\Delta^{q+1}_R$ is contained 
in $V$ as well. By the maximum principle one has for $z\in \Delta^{q+1}_R$  and $s \in B(s_0,\epsilon)$ that
\[
||G^{-1}(g_{s}(z))|| \leqslant \max_{z \in\partial \Delta^q_R \times \partial \Delta_R}{||G^{-1}(g_{s}(z))|| }
\le 
\]
\begin{equation}
\eqqno(max-pr)
 \leqslant 
C\max_{z \in\partial \Delta^q_R \times \partial \Delta_R}{||G^{-1}(g_{*}(z))||_{W^{k,2}(B(s_0,\epsilon),
\mathbb{C}^{m+q+1})} }.
\end{equation}
Here $C$ is the constant from the Sobolev Imbedding Theorem. The mapping 
\[
z \mapsto  G^{-1}(z, F(z,*))
\]
is continuous in a neighborhood of $\partial \Delta^q_R \times
\partial\Delta_R \subset H_q^1(r)$. Set 
\[
\psi_{(z)} (s)= \pi_2 \circ G^{-1}(z,F(z,s)) \in W^{k,2}(B(s_0,\epsilon),\mathbb{C}^m),
\]
where $\pi_2:\Delta^{q+1}_R \times \Delta^m_R \to \Delta^m_R$ 
is the natural projection and consider its Fourier transform $\hat{\psi}$. From the  inequality \eqqref(max-pr) we
see that there exists a constant, namely
\[
M_{s_0}= \max_{z \in \partial \Delta^q_R 
\times \partial \Delta_R}{||{(1+||*||)^k\hat{\psi}_{(z)}(*)}||_{L^2(B(s_0,\epsilon),\mathbb{C}^m)}},
\]
such that for all $z \in \Delta^{q+1}$ one has
\begin{equation}
\eqqno(max-pr1)
||{(1+||*||)^k\hat{\psi}_{(z)}(*)}||_{L^2(B(s_0,\epsilon),\mathbb{C}^m)} \leqslant M_{s_0}.
\end{equation}

Since $S$ is compact one can cover it by a finite number of balls $\{ B(s_0,\epsilon)\}_{s_0 \in J}$ and by 
taking the maximum $M=\max_{s_0 \in J}{M_{s_0}}$ one obtains the inequality
\begin{equation}
||{(1+||*||)^k\hat{\psi}_{(z)}(*)}||_{L^2(S,\mathbb{C}^m)} \leqslant M.
\end{equation}

\smallskip Therefore for all $z \in \Delta^{q+1}_R$ the mapping  $ \psi_{(z)}$ is in $W^{k,2}(S,\mathbb{C}^{n})$  
and consequently the map $ F(z,*)$ is in $W^{k,2}(S,X)$. 
Now let us see that the correspondence $z \mapsto F(z,*)$ is continuous in the variable $z$ 
in Sobolev topology. Indeed, the map $z \mapsto (1+||*||)^k\hat{\psi}_{(z)}(*) $ is a holomorphic 
Hilbert space valued mapping that satisfies the maximum modulus principle, \ie in particular it will continuously 
depend on $z$. The theorem is proved.

\smallskip\qed

\smallskip The following statement gives us one more example of pairs of open sets $\calu\subsetneqq \hat \calu$ 
in a Hilbert 
manifold such that holomorphic mappings extend from $\calu $ to $\hat\calu$; the previous one was $H_q^{\infty}(r)
\subsetneqq \Delta^q\times B^{\infty}$ of Theorem \ref{q1-qinfty}. It shows that  $\hat\calu\deff W^{k,2}(S, \Delta^q
\times \Delta^n)$ is the ``envelope of holomorphy'' of $\calu\deff W^{k,2}(S,H_q^n(r))$. Here $H_q^n(r)\deff (\Delta^q
\times \Delta^n(r))\cup (A^q_{1-r,1}\times \Delta^n)$ stands for the $q$-concave Hartogs figure in $\cc^{q+n}$.

\begin{thm}
\label{env-hol}
Let $\calx$ be a $q$-Hartogs Hilbert manifold. Then every
holomorphic map $F:W^{k,2}(S,H_q^n(r))\to \calx$ extends to a
holomorphic map $\tilde F : W^{k,2}(S, \Delta^q\times \Delta^n) \to \calx$.
\end{thm}
\proof Let $f=(f^q,f^n):S\to \Delta^q\times \Delta^n$ be an element of
$W^{k,2}(S, \Delta^q\times \Delta^n)$. Consider the  analytic
$q$-disk in $W^{k,2}(S,\Delta^q\times \Delta^n)$ defined by 
\begin{equation}
\eqqno(an-disk)
\phi : (z,s)\in \Delta^q\times S \to  \left(h_{f^q(s)}(z),f^n(s)\right),
\end{equation}
where $h_a$ is the following automorphism of $\Delta^q$ interchanging 
$a=(a_1,...,a_q)$ and $0$
\[
 h_a(z) = \left(\frac{a_1-z_1}{1-\bar a_1z_1},..., \frac{a_q-z_q}{1-\bar a_qz_q}\right).
\]
Clearly $\Phi = \phi (\bar \Delta^q)$ is  an analytic disk in $W^{k,2}(S ,
\Delta^q\times \Delta^n)$ possessing the following properties:
\begin{itemize}
\item $\phi(0,s)$ is our loop $f$.

\item For $z\in\d \Delta^q$ one has that $\phi (z, \cdot)(S)\subset A^q_{1-r, 1}
\times \Delta^n$, therefore $\d \Phi\subset W^{k,2}(S,H_q^n(r))$.
\end{itemize}

Consider the $(q+1)$-disk in $W^{k,2}(S,\Delta^q\times \Delta^n)$ defined by
\[
\phi_t(z,s) \deff \phi (z,t,s) =  \left(h_{f^q(s)}(z),
tf^n(s)\right), \qquad |t|<1.
\]
Then
\begin{itemize}
\item $\phi_0 \subset \Delta^q\times \{0\}$.

\item $\phi_1 = \phi $.

\item For all $t\in \Delta$ one has that $\d\Phi_t \subset W^{k,2}(S,
H_q^n(r))$.
\end{itemize}
Therefore for $\eps >0$ small enough the Hartogs figure
\[
H_{q}^1(\eps)\deff \{||z||<1, |t|<\eps \quad\text{or}\quad
1-\eps<||z||<1+\eps, |t|<1\}
\]
is mapped by $\phi$ to
$W^{k,2}(S,H_q^n(r))$.

\smallskip We can now thicken it to an infinite dimensional Hartogs
figure  $H_q^{\infty}(\eps)$ by multiplying it with
$B^{\infty}(\eps)$, where $B^{\infty}(\eps)$ is a ball in
$W^{k,2}(S,\cc^{q+n})$. Decreasing $\eps>0$  if necessary, we
can achieve that the map
\[
\tilde \phi :(z,s,f_2)\to \phi (z,t,s) + f_2(s)
\]
will send $H_q^{\infty}(\eps)$ to $W^{k,2}(S,H_q^n(r))$. Applying
Theorem \ref{q1-qinfty} we extend $\tilde \phi$ to $\Delta^n\times
B^{\infty}$ and therefore $F$ is extended to a neighborhood of
$\phi$. Finally, since $H_q^n(r))$ and $\Delta^q\times \Delta^n$ are
contractible the manifold $W^{k,2}(S, \Delta^q\times \Delta^n)$ is simply
connected for any $S$. This ensures that our extensions give a
single valued holomorphic extension of $F$.

\smallskip\qed

\begin{corol} If $X$ is a $q$-Hartogs complex manifold then every holomorphic
mapping $F:W^{k,2} (S,H_q^n(r))\to W^{k,2} (S,X)$ extends to a holomorphic
mapping $\hat F:W^{k,2} (S,\Delta^q\times \Delta^n)\to W^{k,2} (S,X)$.
\end{corol}
This readily follows from Theorem \ref{env-hol} applied to $q$-Hartogs by 
Theorem \ref{loop-hart} Hilbert manifold $\calx = W^{k,2}(S,X)$.

\newprg[LOOP.almost-h]{Loop spaces of compact complex manifolds are ``almost Hartogs''}
In this subsection we shall explain how to determine Hartogs manifolds among compact complex manifolds
and their loop spaces. Let us start from dimension one. 

\begin{prop}
\label{1-dim}
Every connected Riemann surface except for the Riemann sphere is $1$-Hartogs.
\end{prop}
\proof Let $X$ be a connected Riemann surface. If $X$ is not compact then it is Stein, 
see \cite{Na}. But every Stein manifold $X$ is $1$-Hartogs for the following reason.
Since $X$ can be realized as a closed submanifold in $\cc^N$ every holomorphic
map $F:H_1^1(r)\to X\subset \cc^N$ is given by holomorphic functions which are
extendable to $\Delta^2$ by the classical Hartogs theorem. This gives the extension
of $F$.

\smallskip Now let $X$ be compact. If $X\not=\pp^1$ then its universal cover is either
the disk $\Delta$ or the complex line $\cc$. Since $H_1^1(r)$ and $\Delta^2$ are both
simply connected the extension of $F:H_1^1(r)\to X$ is equivalent to the extension of 
its lift $\tilde F:H_1^1(r)\to \tilde X\equiv \Delta \text{ or } \cc$. And the latter
is just a holomorphic function.

\smallskip Finally, if $X=\pp^1$ then $\cc^2\setminus \{0\}\ni (z_1,z_2) \to [z_1:z_2]$ 
gives an example of a non-extendable map.

\smallskip\qed

Therefore Theorem \ref{loop-hart} implies the following
\begin{corol}
Let $X$ be a connected Riemann surface different from the Riemann sphere. Then the loop space 
$\calw_{\sss^1,X}^{1,2}$ is a $1$-Hartogs Hilbert manifold.
\end{corol}

In \cite{Iv2} we introduced the class $\calg_q$ of $q$-disk convex complex manifolds having
a strictly positive
$dd^c$-closed $(q,q)$-form. The sequence $\{\calg_q\}_{q=1}^{\infty}$ is rather exhaustive: 
$\calg_q$ contains all compact manifolds of dimension $q+1$, see subsection 1.5 in \cite{Iv2} 
for more details.

\smallskip We think that the following statement should be true:

\smallskip\noindent{\slsf Conjecture.} {\it Let $X$ be a complex manifold from the class $\calg_q$. 
Then

\smallskip\sli either $X$ contains a $(q+1)$-dimensional spherical
shell, or $X$ contains an uniruled 

\quad compact subvariety of dimension $q$;

\smallskip\slii or, $X$ is $q$-Hartogs.
}

\begin{rema}
{\bf a)} \rm This was proved in \cite{Iv2} for $q=1$ (in fact this
particular statement was proved already in \cite{Iv1}), and in
\cite{IS} for $q=2$. In the latter paper we proved almost the
assertion stated above (for all values of $q$), but for holomorphic mappings with
zero-dimensional fibers, see Proposition 12 there. The main difficulty in proving
the general case of this conjecture is that the product $\Delta^q\times X$ doesn't
admit a positive $dd^c$-closed $(q,q)$-form in general if even $X$ does. This 
produces difficulties in estimating the areas of graphs of holomorphic maps from
$\Delta^q$ to $X$. These estimates play crucial role in the proofs in \cite{Iv1,Iv2}
and \cite{IS}.

\smallskip\noindent\bf 2. \rm For $q=1$ the item (\sli means just that $X$
contains a shell or a rational curve. For $q=2$ we need to add few explanations
to \cite{IS}. We proved there that a meromorphic map from $H_2^1(r)$
to such $X$ meromorphically extends to $\Delta^3\setminus S$, where $S$
is a complete pluripolar set of Hausdoff dimension zero. If
$S\not=\emptyset$ then $X$ contains a spherical shell of dimension
$3$. Otherwise $S$ is empty. If our map $f$ was in addition
holomorphic on $H_2^1(r)$ then the set $I_{\tilde f}$ of points of
indeterminacy of the extension $\tilde f$ can only be discrete and
then it is clear that for every $a\in I_{\tilde f}$ its full image
$\tilde f[a]$ contains an uniruled analytic set  of dimension
two.

\smallskip\noindent\bf 3. \rm {\slsf Proof of Corollary \ref{sept-0} from the 
Introduction.}
If $X$ is not as in (\sli or (\slii then $X$ is $1$-Hartogs if $\dim X=2$ (resp.
$2$-Hartogs if $\dim X=3$). Theorem \ref{loop-hart} gives that the corresponding
loop spaces are Hilbert-Hartogs as claimed.

\smallskip\noindent\bf 4. \rm Let us make few remarks about a possible interest
in Corollary \ref{sept-0}. The main problem in classifying surfaces of Kodaira class
$VII_0$ with second Betti number $>0$ is to prove that they contain rational
curves. Concerning spherical shells, if such $X$ contains such a shell it is
a deformation of a blown-up Hopf surface. Therefore the $2$-dimensional part
of Corollary \ref{sept-0} specifies that if $X$ has no rational curves and no
shells (\ie is ``totally unknown'') then one can say a good thing about it:
it is $1$-Hartogs itself as well as all its loop spaces. An analogous thing
can be said about a complex structure on $\sss^6$ (if it exists).

\end{rema}

\begin{rema} \rm
\label{chazal}
Finally let us discuss compact manifolds of dimension $q$. Note that they belong to 
$\calg_q$. Indeed, take any volume form on such $X$. It has bidegree $(q,q)$ and is vacously
$dd^c$-closed. It was proved by F. Chazal in \cite{Ch} that every compact manifolds of 
dimension $q$ is {\slsf meromorphically} $q$-Hartogs in the sence that meromorphic 
maps with values in such $X$ extend from $H_q^1(r)$ to $\Delta^{q+1}$. Suppose 
there exists a holomorphic map $f:H_q^1(r)\to X$ which extend to $\Delta^{q+1}$ only
meromorphically. By the same reasoning as in the previous remark we get a uniruled
subpace of $X$ of dimension $q-1$. On the other hand if every $f$ extends holomorphically
then this $X$ is $q$-Hartogs and so are its generalized loop spaces.
\end{rema}

\ifx\undefined\bysame
\newcommand{\bysame}{\leavevmode\hbox to3em{\hrulefill}\,}
\fi

\def\entry#1#2#3#4\par{\bibitem[#1]{#1}
{\textsc{#2 }}{\sl{#3} }#4\par\vskip2pt}

\end{document}